\documentclass[12pt]{article}
\usepackage{amsmath}
\usepackage{amssymb}
\usepackage{epsfig}
\newtheorem{dfn}{Definition}[section]
\newtheorem{prop}[dfn]{Proposition}
\newtheorem{lem}[dfn]{Lemma}
\newtheorem{theo}[dfn]{Theorem}
\newtheorem{coro}[dfn]{Corollary}

\newtheorem{exam}[dfn]{Example}
\newtheorem{rem}[dfn]{Remark}

\newcommand{\bA}{{\mathbb{A}}}
\newcommand{\bC}{{\mathbb{C}}}
\newcommand{\bF}{{\mathbb{F}}}
\newcommand{\bG}{{\mathbb{G}}}

\newcommand{\bP}{{\mathbb{P}}}
\newcommand{\bQ}{{\mathbb{Q}}}

\newcommand{\R}{{\mathbb{R}}}

\newcommand{\bZ}{{\mathbb{Z}}}

\newcommand{\cO}{{\mathcal{O}}}

\newcommand{\cT}{{\mathcal{T}}}

\def\L{{\Lambda}}

\def\a{{\alpha}}
\newcommand{\m}{{\rm m}}
\newcommand{\ra}{{\rightarrow}}
\newcommand{\cN}{{\mathcal N}}

\newcommand{\Pic}{{\rm Pic}}
\newcommand{\Sym}{{\rm Sym}}
\newcommand{\Aut}{{\rm Aut}}
\newcommand{\GL}{{\rm GL}}
\newcommand{\PGL}{{\rm PGL}}
\newcommand{\Spec}{{\rm Spec}}
\newcommand{\End}{{\rm End}}

\author{Brendan Hassett and Yuri Tschinkel}
\title{{\Large Geometry of equivariant 
compactifications of ${\bG}_a^n$}}

\date{February 1999}

\begin{document}

\maketitle

\section{Introduction}

In this paper we begin a systematic study of 
{\em equivariant} compactifications of 
${\bG}_a^n$.  
The question of classifying non-equivariant compactifications
was raised by F. Hirzebruch (\cite{hirzebruch}) and 
has attracted considerable attention since (see  
\cite{furushima}, \cite{peternell}, \cite{mueller-stach}
and the references therein). 
While there are classification results for surfaces and 
non-singular threefolds with small Picard groups,
the general perception is that a complete 
classification is out of reach.

On the other hand, there is a rich theory of equivariant 
compactifications of reductive groups. The classification
of normal equivariant compactifications of 
reductive groups is combinatorial. Essentially, the whole geometry
of the compactification can be understood in terms of (colored) fans. 
In particular, these varieties do not admit moduli. 
For more details see \cite{oda}, \cite{fulton}, \cite{brion} and the
references therein.

Our goal is  to understand {\em equivariant} 
compactifications of $\bG_a^n$. The first step in our approach
is to classify possible $\bG_a^n$-structures on simple
varieties, like projective spaces or Hirzebruch surfaces. 
Then we realize general smooth $\bG_a^n$-varieties as 
appropriate (i.e., equivariant) blow-ups of simple varieties.
This gives us a geometric description of the moduli space of 
equivariant compactifications of $\bG_a^n$.

In section 2 we discuss general properties of
equivariant compactifications of $\bG_a^n$ ($\bG_a^n$-varieties).
In section 3 we classify all possible 
$\bG_a^n$-structures on projective 
spaces $\bP^n$.  
In section 4 we study curves, paying particular attention 
to non-normal examples. In section 5 we carry out our program 
completely for surfaces. 
In particular, we classify all possible $\bG_a^2$-structures
on minimal rational surfaces.  In section 6 we turn to threefolds.
We give a classification of $\bG_a^3$-structures on 
smooth projective threefolds with Picard group of rank 1. 
Each section ends with a list of examples and open questions.

\

{\bf Acknowledgements.} We are grateful to V. Batyrev, 
S. Bloch, T. Defernex, V. Ginzburg, D. Kazhdan and
B. Totaro for useful discussions. The first author
was partially supported by an NSF postdoctoral fellowship. 
The second author
was partially supported by the NSA.

\section{Generalities}

We work in the category of algebraic varieties over 
$F=\overline{\bQ}$.

\subsection{Definitions}

\begin{dfn}
Let ${\bG}$ be a connected linear algebraic group.
An algebraic variety $X$ admits a (left) $\bG$-action if
there exists a morphism 
$\varphi \,:\, \bG \times X\ra X$, satisfying the standard
compatibility conditions.  
A ${\bG}$-variety $X$  is a variety  with a {\em fixed}
(left) $\bG$-action such that 
the stabilizer of a generic point is trivial and the 
orbit of a generic point is dense.
\end{dfn}
For example, a normal $\bG_m^n$-variety is a toric variety. 
\begin{dfn}
\label{dfn:equiv}
A morphism of $\bG$-varieties is a morphism 
of algebraic varieties commuting with the $\bG$-action.
A $\bG$-{\em isomorphism} is an isomorphism in the category of
$\bG$-varieties. A $\bG$-equivalence is 
a diagram 
$$
\begin{array}{ccc}
{\bG}\times X_1 & \stackrel{(\alpha,j)}{\longrightarrow}      & {\bG}\times X_2 \\
\downarrow      &                               &  \downarrow\\
X_1             &             \stackrel{j}{\longrightarrow} &   X_2 
\end{array}
$$
where $\a \in \Aut (\bG)$ and $j$ is an isomorphism 
(in the category of algebraic varieties).      
\end{dfn}
Clearly, every $\bG$-isomorphism is a $\bG$-equivalence.
We shall omit $\bG$ if the group is understood. 

One of our main observations in this paper is that
classification of simple $\bG$-varieties up to 
equivalence, even projective
spaces,  is a non-trivial problem. 
This is in marked contrast to the situation for toric varieties. 
A toric variety admits a unique structure as a  $\bG_m^n$-variety
(up to equivalence). Here is a sketch of the argument for
a projective toric variety  $X$: 
Consider the connected component
of the identity $\Aut(X)^0$ of
the automorphism group $\Aut(X)$. 
This is an algebraic group which acts trivially on the 
Picard group  of $X$. We pick a very ample line bundle
and consider $\Aut(X)^0$ as a closed subgroup of the corresponding 
group $\PGL_N$. In particular, $\Aut(X)^0$ is a linear algebraic
group. 
The key ingredient now is the statement that 
all maximal tori in $\Aut(X)^0$ are conjugate
(cf. 11.4 in \cite{borel-lin-alg}).  
Evidently, the maximal torus
acts faithfully on $X$, and the action 
has a dense open orbit. This proves the claim.

Every ${\bG}$-variety contains an open subset
isomorphic to ${\bG}$. 
We denote by $D$ the complement of this open subset and call
$D$ the {\em boundary}. 
If $X$ is normal, Hartog's theorem implies that 
$D$ must be a divisor (the complement to $D$ is affine). 
Otherwise, we normalize and observe that 
the normalization is an isomorphism over $\bG$.

\subsection{Line bundles and linearizations}

\begin{prop}\label{prop:linearization}
Let $X$ be proper and normal 
algebraic variety. 
Then the action of ${\bG}_a^n$ on the Picard group $\Pic(X)$ 
is trivial and every line bundle on $X$
admits a unique linearization, up to scalar multiplication.  
\end{prop}

{\em Proof.} Uniqueness 
follows from Prop. 1.4 p. 33 in \cite{mumford}.
(The only relevant hypothesis is 
that $X$ is geometrically reduced,
which is part of our assumptions.) 
The proof of Prop. 1.5 p. 34 implies that 
every line bundle on $X$ has a linearization. 
(Here we use that $X$ is proper and 
that the Picard group of the group $\bG$ is trivial.
To get the fact that the action of $\bG$ on the $\Pic(X)$ is
trivial we need the normality of $X$.)

\begin{coro}\label{coro:linear} 
Retain the assumptions of Proposition \ref{prop:linearization}.
Consider a basepoint free
linear series $W\subset H^0(X,L)$, 
stable under the action of $\bG_a^n$.   
Then the induced map $f\,:\, X\ra \bP(W^*)$ (here we use
the {\em geometric} convention) is $\bG_a^n$-equivariant. 
\end{coro}
 
\begin{theo}\label{theo:pic}
Let $X$ be a proper and normal $\bG_a^n$-variety.
Then 
$\Pic(X)$ is freely generated by the classes
of the irreducible components $D_j$ ($j=1,...,t$) 
of the boundary divisor $D$. 
The cone of effective Cartier divisors 
$\L_{\rm eff}(X)\in \Pic(X)_{\R}$ is 
given by
$$
\L_{\rm eff}(X)= \oplus_{j=1}^t \R_{\ge 0} [D_j].
$$
\end{theo}

{\em Proof.}
Choose an effective divisor $A\subset X$ and 
consider the representation of $\bG_a^n$ on 
the projectivization $H^0(X, {\cO}(A))$
(here we use that $X$ is normal). 
This representation has a fixed point corresponding 
to an effective divisor supported at the boundary. 
To show that there are no relations between the classes
$[D_1],...,[D_t]$ it suffices to observe that there exist
no functions without zeros and poles on $\bG_a^n$.

\begin{rem}
Every effective cycle on a $\bG$-variety 
is rationally equivalent to a cycle supported on the boundary. 
(Here we are using the fact the $\bG$ is affine.)
\end{rem}

\subsection{Vector fields and the anticanonical line bundle}

\begin{theo} \label{theo:anticanonical}
Let $X$ be a smooth and proper $\bG$-variety.
Then the anti-canonical class is a sum of classes
of the irreducible components of the boundary $D$ 
with coefficients which are all $\ge 1$. 
If $\bG={\bG}^n_a$ then the coefficients are all $\ge 2$.
\end{theo}

{\em Proof.}
We first introduce some general terminology and exact sequences.
Let $X$ be a smooth variety and $D$ a normal-crossings divisor.
Then there are two natural exact sequences
\begin{eqnarray*}
0 & \rightarrow & {\cT}_X\left< -D\right> 
\rightarrow {\cT}_X \rightarrow
{\cN}_{D/X} \rightarrow 0 \\
0 & \rightarrow & {\cT}_X(-D) \rightarrow 
{\cT}_X\left< -D \right> \rightarrow
{\cT}_D \rightarrow 0 
\end{eqnarray*}
where ${\cT}_X\left< -D \right>$ 
denotes the {\em vector fields with logarithmic
zeros}.  If $D$ is given locally by the equation $t=0$ these
take the form $t\frac{\partial}{\partial t}, 
\frac{\partial}{\partial s}$.

Let $v\in H^0(X,{\cT}_X)$ be a vector field.  By definition,
$v$ {\em vanishes normally to order one along $D$} 
if its image in $H^0(D,{\cN}_{D/X})$ is zero.
If $v$ arises from the action of a one-parameter group
stabilizing $D$ then it vanishes normally to order one along $D$.  
If $v$ vanishes normally to order one along $D$, we can consider 
the corresonding element $w\in H^0(D,{\cT}_D)$.  If $w=0$ 
then we say that $v$ {\em vanishes to order one}.  
If $v$ arises from a one-parameter group fixing $D$ then
it vanishes to order one.  
Generally, if $v$ vanishes to order $N-1$
(resp. normally to order $N$)
then we can consider the corresponding element
$w\in H^0(D,{\cN}_{D/X}(-(N-1)D))$ (resp.
$H^0(D,{\cT}_D(-(N-1)D))$).  
If $w=0$ then we say that $v$ vanishes normally
(resp. vanishes) to order $N$ along $D$;  in 
particular, $v$ may be regarded as a section of 
$H^0(X,{\cT}_X\left<-D \right>(-(N-1)D))$ (resp.
$H^0(X,{\cT}_X(-ND))$).  

Now assume that $X$ is a smooth proper ${\bG}$-variety
of dimension $n$ and $D_i$ a component of its boundary.  
Let $\{ v_1,\ldots ,v_n \} \in H^0(X,{\cT}_X)$ be 
invariant vector fields spanning the Lie algebra of $\bG$, and
$M_j$ (resp. $N_j$) the order of vanishing (resp. normal
vanishing) of $v_j$.  By definition, $N_j=M_j$ or $M_j+1$.
Consider the exterior power 
$\sigma=v_1 \wedge v_2\wedge 
\ldots \wedge  v_n \in H^0(X,\Lambda^n{\cT}_X);$
we bound (from below) the order of 
vanishing $a_i$ of $\sigma$ along $D_i$.  
Evidently $a_i\ge M_1 + \ldots +M_n$ 
but in fact slightly more is true.
Using the adjunction isomorphism
$$\Lambda^n {\cT}_X|{D_i}=\Lambda^{n-1} {\cT}_{D_i} 
\otimes {\cN}_{D_i/X}$$
we obtain
$$a_i\ge \min_{j=1,\ldots,n}(M_1 + \ldots +M_{j-1}+N_j+M_{j+1}+ \ldots + M_n).$$
For instance, since each $N_j>0$ we obtain that $a_i>0$.
It follows that the canonical divisor is linearly equivalent to 
$-\sum_i a_i D_i$, where each $a_i>0$.  

Now assume that $X$ is a smooth proper ${\bG}^n_a$-variety
and $D_i$ is a component of its boundary.  Let $v_1$ be a vector
field arising from the action of a one-parameter subgroup
that fixes $D_i$, so that $v_1$ vanishes  to order
one and $M_1\ge 1$.  We claim $v_1$ necessarily vanishes normally
to order {\em two} i.e. $N_1\ge 2$.  Consider the resulting element
$w \in H^0(D_i,{\cN}_{D_i/X}(-D_i))=\End({\cN}_{D_i/X})$, 
which exponentiates to give the induced ${\bG}^1_a$ action
on the normal bundle.  Since ${\bG}^1_a$ has no characters, this
action is trivial and $w=0$.  The previous inequality guarantees
that $a_i\ge 2$.  Hence the canonical divisor is linearly equivalent
to $-\sum_i a_i D_i$ where each $a_i>1$.

\begin{coro}
\label{coro:Pncrit}
Let $X$  be a smooth projective ${\bG}^n_a$-variety
with irreducible boundary $D$, which is fixed under the action.
Then $X=\bP^n$.
\end{coro}

{\em Proof.}
We see that $M_j \ge 1$ for each $j$.
It follows that $K_X=-rD$ where $r\ge n+1$, i.e. $X$ is Fano
of index $r\ge n+1$.  Hence $X={\bP}^n$ and $r=n+1$.

\begin{coro}
\label{coro:Qncrit}
Let $X$  be a smooth projective ${\bG}^n_a$-variety
with irreducible boundary $D$.  Assume that the
subgroup of $\bG_a^n$ fixing  $D$ has dimension $n-1$.  
Then $X={\bP}^n$ or $X=Q_n$, the quadric 
hypersurface of dimension $n$.  
\end{coro}

{\em Proof.}
After reordering, we obtain that $M_j \ge 1$ for $j=1,\ldots,n-1$.  
It follows from the proof of Theorem \ref{theo:anticanonical} that
$N_j\ge 2$ for $j=1,\ldots,n-1$ and $N_n\ge 1$.  In particular, 
$r \ge n$, i.e., $X$ is Fano of index $r\ge n$.  Hence $X={\bP}^n$
or $Q_n$.

\subsection{A dictionary}

Let $L({\bG}_a^n)$ and $U({\bG}_a^n)$ denote the Lie algebra
and the enveloping algebra of ${\bG}_a^n$.  Since ${\bG}_a^n$
is commutative $U({\bG}_a^n)$ is isomorphic to a
polynomial ring in $n$ variables.  Let $R=U({\bG}_a^n)/I$ and
assume that $\Spec(R)$ is supported at the origin;
we use $\ell(R)$ and $m_R$ to denote the length and the
maximal ideal of $R$.  Since
$I$ contains all the homogeneous polynomials 
of sufficiently large degree $d$, 
the elements of $L({\bG}_a^n)$ act (via the regular
representation) as nilpotent matrices on $R$.  
We can exponentiate to obtain an algebraic representation
$$ \rho_R: {\bG}_a^n \rightarrow \Aut_F(R)$$
of dimension $\ell(R)$.

For concreteness, we introduce some additional notation.
We consider 
$L({\bG}_a^n)$ as a vector space over $F$ with a distinguished 
basis $S_j=\frac{\partial}{\partial x_j}$
so that $U({\bG}_a^n)=F[S_1,\ldots,S_n]$.  
The representation $\rho_R$ is obtained 
by multiplying by
$\exp (x_1S_1+\ldots + x_nS_n) \in R.$
We shall use the following basic fact:
\begin{lem}
There is an $F$- basis $\{\mu_1=1,\mu_2,\ldots,\mu_{\ell(R)} \}$
for $R$ such that
each $\mu_j$ (with $j>1$) is a monomial of the form
$S_i\mu_k$ for some $i,k$.
\end{lem}
{\em Proof}.  This is proved using Gr\"obner basis techniques 
\cite{Eisenbud} chapter 15.  With respect to homogeneous
lexicographic order, the initial ideal $\mathrm{Init}(I)$ is monomial.  
The monomials not contained in $\mathrm{Init}(I)$ form a basis for $R$ 
(cf. proof of \cite{Eisenbud} Theorem 15.17.)  These satisfy
the desired properties.

The lemma implies that we may expand
$$\exp(x_1S_1+\ldots +x_nS_n)=
\sum_{i=1}^{\ell(R)} f_j(x_1,\ldots,x_n)\mu_j;$$
note that $f_1=1$.  
\begin{prop}
The coordinate functions satisfy the following:
\begin{enumerate}
\item{ If $\mu_j=S_i\mu_k$ then $f_k=\frac{\partial f_j}{\partial x_i}.$}
\item{The polynomials 
$f_1,\ldots,f_{\ell(R)}$ form a basis for the space of  solutions
of the system of equations 
$$g[f(x_1,\ldots,x_n)]=0 \text{ for each }g\in I.$$}
\end{enumerate}
\end{prop}
{\em Proof.} We interpret the equation
$$\frac{\partial}{\partial x_i}\exp(x_1S_1+\ldots +x_nS_n)=
S_i\exp(x_1S_1+\ldots +x_nS_n)$$
in terms of our expansion for the exponential.
This implies the first assertion and that the $f_j$ are solutions
to the system of equations.

We next show that the $f_j$ are independent.   By the first
assertion, each $f_j$ contains a monomial 
$x_1^{m(1)} \ldots x_n^{m(n)}$ where 
$\mu_j=S_1^{m(1)}\ldots S_n^{m(n)}$.  Homogeneous lexicographic
order on the monomials in the $S_i$ induces an order on the 
monomials in the $x_i$;  simply regard $x_i$ as $S_i^{-1}$.  
With respect to this order, $x_1^{m(1)}\ldots x_n^{m(n)}$ is clearly 
the leading term of $f_j$. 
Since their leading terms are distinct 
the $f_j$ are linearly independent.     

To complete the proof we show that the solution space 
$$V:=\{f(x_1,\ldots,x_n):g[f]=0 \text{ for each }g\in I\}$$
has dimension $\ell(R)$.  
Since $I$ contains all the polynomials of degree $d$, 
each solution is polynomial with total degree $<d$.  
There is a natural pairing between $F[x_1,\ldots,x_n]$ and
$F[\frac{\partial}{\partial x_1},\ldots,\frac{\partial }{\partial x_n}]$
$$\left<g,f \right>=g[f]|_{(0,\ldots,0)}.$$
This induces a perfect pairing between homogeneous polynomials and
operators of a given degree.  Note that
$V= \{ f: \left< g,f\right>=0 
\text{ for each } g\in I\}$
which implies that $\dim V = \ell(R)$.

We collect some basic properties of $\rho_R$:
\begin{prop}
The ${\bG}_a^n$-representations $\rho_R$ and $V$ are dual,
and $\rho_R$ has a nondegenerate orbit 
(i.e., a cyclic vector) in $R$.
The representation 
$\rho_R$ is faithful iff $\Spec(R) \subset \Spec(F[S_1,\ldots,S_n])$
is nondegenerate.  
\end{prop}
The vector space $V$ has a natural
$\bG_a^n$-action by translations. 
The first statement is clear from the preceding discussion.  
Consider the orbit of $1\in R$
$$\rho_R(x_1,\ldots,x_n)\cdot 1
=\sum_{i=1}^{\ell(R)} f_j\mu_j.$$
The linear independence of the $f_j$ implies that this  orbit 
is nondegenerate in $R$.
The final statement is clear;  indeed, $S_i$ acts trivially iff 
$\Spec(R)\subset \{ S_i=0 \}$. 

\begin{rem}
Note that $V=\rho^*_R$ has a natural structure as an $R$-module
of length $\ell(R)$.  Indeed, it coincides with the dualizing
module $\omega_R$.  This can be seen using Macaulay's method of 
inverse systems (see \cite{Eisenbud} chapter 21.2).  It follows that
$\rho_R$ is self-dual iff $R$ is Gorenstein.
\end{rem}  
 
A translation invariant
subspace $V\subset F[x_1,\ldots,x_n]$ of dimension $\ell$
corresponds to a
representation of ${\bG}^n_a$ on $F^{\ell}$ with a fixed cyclic vector.
Consider the ideal $I$ of constant coefficient differential operators
annihilating $V$.  Setting 
$$R= F[\frac{\partial}{\partial x_1}, \ldots,
\frac{\partial}{\partial x_n}]/I$$
we obtain a scheme $\Spec(R)$ of length $\ell$ supported
at the origin.  

We obtain a dictionary:
\begin{theo}
There is a one-to-one correspondence among the following:
\begin{enumerate}
\item{
subschemes $\Spec(R) \subset \Spec( 
F[\frac{\partial}{\partial x_1},\ldots,
\frac{\partial }{\partial x_n}])$
supported at the origin of length $\ell=\ell(R)$;}
\item{
translation invariant
subspaces $V \subset F[x_1,\ldots ,x_n]$ of dimension $\ell$;}
\item{isomorphism classes of pairs $(\rho,v)$ such that 
$\rho:{\bG}_a^n \rightarrow \GL_{\ell}$
is a representation and $v$ is a cyclic vector (i.e.,
$\rho({\bG}^n_a)\cdot v$
is nondegenerate).}
\end{enumerate}
\end{theo}

We now turn to a case of particular interest.  
Assume that $S_1,\ldots,S_n$ form a basis 
for the maximal ideal $m_R$, i.e.,
$\ell(R)=n+1$.  Then the corresponding
representation $\rho_R$ is faithful with {\em dense} orbit on
${\bP}(R)$.  More intrinsically, for any Artinian local $F$-algebra,
exponentiating the action of $m_R$ on $R$ yields a representation
$\rho_R$ of ${\bG}_a^{\ell(R)-1}$ with dense open orbit on 
${\bP}(R)={\bP}^{\ell(R)-1}$.  Conversely, any faithful 
representation $\rho$ with dense projective orbit yields
an Artinian local $F$-algebra of length $\ell$, i.e.,
the differential operators with constant coefficients
annihilating its coordinate functions.  
We summarize this discussion in the proposition:
\begin{prop}
\label{prop:projclass}
The following are equivalent:
\begin{enumerate}
\item{Artinian local $F$-algebra $R$ of length $\ell=\ell(R)$, up to 
isomorphism;}
\item{equivalence classes of ${\bG}_a^{\ell-1}$-structures on
${\bP}^{\ell-1}$.}
\end{enumerate}
\end{prop}

\subsection{Examples and questions}

\begin{enumerate}
\item{Not every point in the boundary $D$ is contained
in the closure of a 1-parameter subgroup.
Construction: Blow up $\bP^2$ in a point at the boundary. Blow up
again a point in the exceptional divisor. 
Every 1-parameter subgroup in $\bP^2$ is a line.} 
\item{
The theorem \ref{theo:anticanonical} fails for 
non-equivariant compactifications of
$\bG_a^n$. For example, let $X\subset \bP^2\times \bP^2$
be a hypersurface of bidegree $(1,d)$ with $d\ge 4$. 
Then the anticanonical class ${\cal O}(2,3-d)$ is not
contained in the interior of the effective cone.} 
\item Suppose $X$ is a smooth projective $\bG_a^n$-variety
with finitely many $\bG_a^n$-orbits. Is  $X$ rigid as an algebraic
variety?
\end{enumerate}

\section{Projective spaces}
In this section we study ${\bG}^n_a$-structures
on projective spaces. Notice that every $\bP^n$ has
a distinguished structure as a $\bG_a^n$-variety. The translation
action on the affine space $\bA^n$ extends to an action on $\bP^n$,
fixing the hyperplane at infinity. 
We denote this action by $\tau_n$. It corresponds to the 
Artinian ring $F[S_1,...,S_n]/[S_jS_j, \,\, i,j =1,...,n] $.
It is easy to see that {\em every} $\bG_a^n$-structure on 
$\bP^n$ admits a specialization to $\tau_n$.

\

In the following propositions, we classify ${\bG}^n_a$-structures
on projective spaces 
of small dimension, {\em up to equivalence} (cf. \ref{dfn:equiv}).
The first natural invariant is the Hilbert-Samuel function of
the corresponding Artinian ring $R$, defined by $\chi_R(k)=
\ell(m_R^k/m_R^{k+1})$.

\begin{prop}
There is a unique ${\bG}^a_1$-structure on ${\bP}^1$.
\end{prop}
{\em Proof.}  It is a consequence of Proposition 
\ref{prop:projclass}.  

\begin{prop}\label{prop:p2}
There are two distinct $\bG_a^2$-structures on $\bP^2$. 
They are given by the following representations of $\bG_a^2$:
$$
\tau_2({a_1,a_2})=\left( \begin{matrix}
1 & 0 & a_2  \\
0 & 1 & a_1 \\
0 & 0 & 1
\end{matrix}\right)$$
and
$$\rho({a_1,a_2})=\left( \begin{matrix}
1 & a_1 & a_2+\frac{1}{2}a_1^2  \\
0 & 1 & a_1 \\
0 & 0 & 1
\end{matrix}\right).$$ 
They correspond to the quotients of $F[S_1,S_2]$ by the ideals
$I_1=[S_1S_2,S_2^2,S_1^2]$ and $I_2=[S_1S_2,S_2-S_1^2]$.
\end{prop}
{\em Proof.}  
It suffices to classify Artinian local $F$-algebras $R$ of length
three up to isomorphism.  
If the tangent space has dimension two then
$R=F[S_1,S_2]/I_1$.   
If the tangent space to $R$ has dimension one
then $R=F[S_1,S_2]/I_2$ (it is clear 
that this is the only one).  
These correspond to the 
representations $\tau_2$ and $\rho$ respectively.

\begin{prop}\label{prop:p3}
There are four distinct $\bG_a^3$-structures on $\bP^3$. 
They correspond to the quotients of $F[S_1,S_2,S_3]$ by
the following ideals:
\begin{eqnarray*}
I_1&=&[S_1^2-S_2,S_1S_2-S_3,S_1S_3]\\
I_2&=&[S_1^2-S_2,S_1S_2,S_1S_3]\\
I_3&=&[S_1^2,S_1S_2-S_3,S^2_2]\\
I_4&=&[S_1^2,S_1S_2,S_2^2,S_2S_3,S_3^2,S_1S_3].
\end{eqnarray*}
\end{prop}
{\em Proof.}  If the tangent space to $R$ has dimension one or three,
then $R$ is necessarily the quotient of the polynomial ring by
$I_1$ or $I_3$.  If the tangent space has dimension two then
$\ell(m^2_R/m^3_R)=1$.  Consider the nonzero symmetric quadratic form
$$
q\,:\, m_R/m^2_R \times m_R/m^2_R \ra m^2_R / m^3_R.$$
Let $S_3$ generate $m^2_R$.  Then there exist $S_1,S_2 \in m_R$
spanning $m_R/m_R^2$ such that the quadratic form equals
$S_1S_2$ or $S_1^2$.

\begin{prop}\label{prop:p4}
There are ten distinct $\bG_a^4$-structures on $\bP^4$. 
They correspond to the quotients of $F[S_1,S_2,S_3,S_4]$ by
the following ideals:
\begin{eqnarray*}
I_1&=&[S_1^2-S_2,S_1S_2-S_3,S_1S_3-S_4,S_2S_3,S_1S_4]\\
I_2&=&[S_1^2-S_3,S_1S_2,S_2^2,S_1S_3-S_4,S_1S_4]\\
I_2&=&[S_1^2-S_3,S_1S_2,S_2^2,S_1S_3-S_4,S_1S_4]\\
I_3&=&[S_1^2-S_3,S_1S_2,S_2^2-S_3,S_1S_3-S_4,S_1S_4]\\
I_4&=&[S_1^2-S_3,S_1S_2-S_3,S_2^2,S_1S_3-S_4,S_1S_4]\\
I_5&=&[S_1^2-S_3,S_1S_2-S_4,S_2^2,S_1S_3,S_1S_4,S_2S_3]\\
I_6&=&[S_1^2-S_3,S_2^2-S_4,S_1S_2,S_1S_3,S_2S_4]\\
I_7&=&[S_1^2-S_4,S_2S_3-S_4,S_1S_2,S_1S_3,S_2^2,S_3^2,S_1S_4]\\
I_8&=&[S_1^2,S_2^2,S_3^2,S_1S_2-S_4,S_1S_3,S_2S_3]\\
I_9&=&[S_1^2-S_4,S_2^2,S_3^2,S_1S_2,S_1S_3,S_1S_4,S_2S_3]\\
I_{10}&=&[S_iS_j,\,\, i,j=1,...,4]
\end{eqnarray*}
\end{prop}
{\em Proof.} 
We consider the possible shapes of the Hilbert-Samuel function
$\chi_R$. In the cases where 
$\chi_R(1)=1$ or $4$ it is clear that the only possibilities
are $I_1$ and $I_{10}$ respectively. 

Assume that $\chi_R(1)=2$, $\chi_R(2)=1$ and $\chi_R(3)=1$. 
Choose $S_4\in m^3_R,S_3\in m^2_R,$ and $S_1,S_2\in m_R$
which span the maximal ideal of the graded ring associated to $R$.  
For a suitable choice of $S_2$ and $S_4$, we may assume that $S_2S_3=0$
and $S_1S_3=S_4$.  Consider the map $s_2:m_R/m_R^2\ra m_R^2/m_R^3$
induced by multiplying by $S_2$.  If $s_2=0$ then we
may choose $S_3$ so that $S_1^2=S_3$;  we obtain $I_2$.  
If $s_2 \ne 0$ and $S_2 \not \in \ker(s_2)$ 
then we may choose $S_1$ so that $S_1S_2=0$.  However, $S_4 \ne 0$
implies $S_1^2 \ne 0$, so after rescaling $S_1,S_2,$ and $S_4$ we 
obtain $I_3$.  If $s_2 \ne 0$ and $S_2 \in \ker(s_2)$ then
(after rescaling $S_3$) we obtain $S_1S_2=S_3$.  Again,
after rescaling $S_1$ and $S_2$,
$S_1^2 \ne 0$ and we obtain $I_4$.   

Assume that $\chi_R(1)=2$ and $\chi_R(2)=2$. 
Then we choose $S_1,...,S_4$ such that $S_3 $ and $S_4$ are
in $m_R^2$ and $S_1,S_2$ are independent modulo $m_R^2$. 
The ring structure on $R$ is determined by the 
vector-valued quadratic form
$$
q\,:\, m_R/m_R^2 \times m_R/m_R^2 \ra m_R^2/m_R^3.
$$
This corresponds to choosing a 
codimension 1 subspace of ${\rm Sym}^2(m_R/m_R^2)$.
Up to changes of coordinates in $S_1$ and $S_2$, 
each such subspace is spanned by vectors 
$\{S_1^2,S_2^2\}$ and $\{S_1S_2,S_1^2\}$.  This gives the 
cases $I_5,I_6$. 

Assume that $\chi_R(1)=3$ and $\chi_R(2)=1$. 
This corresponds to $I_7,I_8,I_9$. The ring structure
is determined by the quadratic form $q$ (with values in $F$), which 
has rank $3,2$ or $1$. 

\begin{prop}\label{prop:p5}
There are finitely many distinct $\bG_a^5$-structures on $\bP^5$. 
\end{prop}

{\em Proof.} These are written out explicitly in
Suprunenko \cite{supr} pp. 136-150.

The above discussion mirrors the classification of
algebras of commutative nilpotent matrices in the book \cite{supr}. 
(Notice a misprint 
in the classification of algebras corresponding to 
$\bG_a^3$-structures on $\bP^3$ on page 134). 
The arguments in this book
yield a classification of Artinian algebras of length $n+1$
with the following Hilbert-Samuel functions (though the 
author does not make this explicit):
$$
\begin{array}{cccccc}
\chi_R(0) & \chi_R(1)   & \chi_R(2)     & \chi_R(3) & ... &\chi_R(n)\\
1          & n          &  0            &  0 & ...  &    0     \\ 
1          & 1          &  1            &  1 & ...  &    1  \\
1          & 2          &  1            &  1 & ...  & 0\\
1          & 3          &  1            &  1 & ...  & 0 \\
1          & 2          &  2            &  1 & ...  & 0 \\
\end{array}
$$
In particular,  there are finitely many algebras with each of these
Hilbert-Samuel functions $\chi_R$. 
This suffices to obtain a complete classification of 
$\bG_a^n$ structures on $\bP^n$ for $n\le 5$.

Beginning with dimension $6$ we obtain moduli. 
As an example, let us consider 
Artinian rings with Hilbert-Samuel
function of the shape $(1,n-k,k, 0 ,...,0)$ (for suitable $k$). 
These correspond to $k$-dimensional spaces of quadratic forms
in $n-k$ variables (up to coordinate transformations
of the $n-k$ variables). The quadratic forms are
obtained by dualizing the natural map
$$
\Sym^2(m_R/m_R^2)\ra m_R^2/m_R^3. 
$$
For example, if $n=6$ and $k=2$ the moduli space is birational
to the moduli space of elliptic curves.
\begin{exam}\label{prop:p6}
There exists at least one $1$-parameter family of inequivalent 
$\bG_a^6$-structures on $\bP^6$. 
\end{exam}

If $n=8$ and $k=3$ we get the moduli space of genus 5 curves. 
If $n=9$ and $k=3$ we obtain the moduli space of K3 surfaces
of degree 8. 
The appearance of these K3 surfaces and the genus 5 curves
is quite interesting. Can it be explained  
geometrically, in terms of birational maps between 
different $\bG_a^n$-structures?

For a good general introduction to Artinian rings 
(and many further references)
see \cite{iarrobino}.

\begin{prop}
The projective space $\bP^n$ admits a unique $\bG_a^n$-structure with
finitely many orbits. It 
corresponds to the Artinian ring
$F[S_1,...,S_n]/I$, where 
\begin{eqnarray*}
I &=&[S_1^2-S_2,S_1S_2-S_3,...,S_1S_{n-1}-S_n,
S_iS_j,\,\, i+j >n] \\
 &\simeq & [S_1^i-S_i, S_iS_j, \,\, i+j>n] . 
\end{eqnarray*}
\end{prop}

{\em Proof.}
There is a unique fixed point under the $\bG_a^n$-action.
Projecting from it
gives a $\bP^{n-1}$ with finitely many $\bG_a^{n-1}$-orbits.
By the  inductive hypothesis, this 
$\bP^{n-1}$ has the indicated structure.
The Artinian ring $R(\bP^{n-1})$  for $\bP^{n-1}$ 
is a quotient of the Artinian ring $R(\bP^{n})$.
Let $S_n\in R(\bP^{n})$ be a non-zero 
element mapped to $0\in R(\bP^{n-1})$. 
Since  $R(\bP^{n-1})\simeq F[T_1]/[T_1^n]$,  
there exists an element $S_1\in R(\bP^{n-1})$ such that
$S_1^{n-1}\neq 0$.  Then $S_1^n=cS_n$ for some constant $c\in F$.
If $c=0$ the action has infinitely many orbits. Otherwise, 
(after rescaling) we obtain the desired ideal. 

\subsection{Examples and questions}

\begin{enumerate}
\item Untwisting different actions on $\bP^2$:
Take $\bP^2$ with the $\rho$-action and
choose a generic one-parameter subgroup.
Let $C$ be the conic obtained as the  
closure of a generic orbit of the one-parameter
subgroup. 
The curve $C$ is tangent to the line at infinity at the fixed
point $p$. Blow up 
the fixed point   {\em on $C$} 3 times. Contract the strict transforms
of the line at infinity, and the first two exceptional curves. 
\item Give an explicit factorization
for $\bG_a^n$-equivariant birational automorphisms
of the projective space $\bP^n$. 
\item Give a dictionary between $\bG_a^n$-structures
on smooth quadrics $Q_n$ and certain Artinian rings (with additional 
structure). 
\item Is there an {\em irreducible} variety parametrizing all
${\bG}_a^n$-structures on ${\bP}^n$? Probably not.   
\end{enumerate}

\section{Curves}

\begin{prop} \label{prop:p1} 
Every {\em smooth} proper $\bG_a^1$-variety
is isomorphic to $\bP^1$ with the standard 
translation action $\tau_1$.
\end{prop}

{\em Proof.} Exercise.

\begin{lem}\label{lem:filt}
Let $V$ be the standard representation of
$\bG_a^1$. Then $\Sym^n(V)$ 
has a filtration 
$$
0\subset {\rm F}_0\subset {\rm F}_1 \subset ... \subset 
{\rm F}_n =\Sym^n(V)
$$
which is compatible with the $\bG_a^1$-action and 
such that
${\rm F}_i\simeq \Sym^i(V) $  and ${\rm F}_{i+1}/{\rm F}_i $
is the 1-dimensional trivial representation. 
Furthermore, every stable subspace of $\Sym^n(V)$ 
arises in this way.
\end{lem}

{\em Proof.} Exercise.

\begin{rem}
The Jordan canonical form gives us a complete description of
representations of $\bG_a^1$. They are isomorphic to direct
sums of the representations $\Sym^n(V)$. 
\end{rem}

\begin{prop}\label{prop:p1-all}
Every proper $\bG_a^1$-variety
$C$ {\em with an equivariant projective embedding}
is isomorphic to $\bP^1$ embedded by a complete
linear series  with the translation action.
\end{prop}

{\em Proof.}
Clearly, the normalization of the 
curve $C$ is isomorphic to $\bP^1$.
Furthermore, the normalization map $\nu \, :\, \bP^1\ra C$ 
is equivariant and an isomorphism away from the fixed point 
$P_{\infty}\in \bP^1$. The morphism $\nu$ is given by some
basepoint free linear series $W$ on $\bP^1$, which is stable
under the $\bG_a^1$-action. 
In particular, 
$W\subset H^0({\cal O}_{\bP^1}(n))\simeq \Sym^n(V)$
(where $V$ is the standard representation).
By the previous Lemma, each stable proper subspace of 
$\Sym^n(V)$ corresponds to a linear series on $\bP^1$
{\em with basepoints}. This concludes the proof.

Our next goal is to classify proper 1-dimensional 
$\bG_a^1$-varieties $C$. Clearly, the normalization 
$\tilde{C}$
has to be isomorphic to $\bP^1$ 
with the standard action $\varphi$
and with the conductor-ideal vanishing at the fixed point
$t=0$. Hence it suffices to classify conductor-ideals   
$I\subset F[t]|_{t=0}$, stable under the group action.

\begin{theo}\label{theo:1-dim-class}
The only conductor-ideals $I$ which are stable under the
group action  $\varphi$ are preimages of some semigroup $\Sigma
\subset \bZ_{<0}$
under the valuation homomorphism. 
\end{theo}

{\em Proof.} Consider the complete 
local ring $F[[t]]$ with maximal
ideal $\m$.  
We make the identification 
$$
\frac{\m}{\m^{n+2}}\simeq \Sym^n(V)
$$
(it follows from the definition of the action $\varphi\,:\, 
t\mapsto t\cdot(1-at+a^2t^2-...)$ extended to the completion). 
By lemma \ref{lem:filt}, all subspaces stable under
the action of  $\bG_a$ coincide with  preimages of subsets of 
the valuation group. 

\subsection{Examples and questions}

\begin{enumerate}
\item 
Let $C\in \bP^2$ be a cuspidal cubic plane curve. 
Then it is a proper $\bG_a^1$-variety
which does not admit
${\bG}_a^1$-equivariant projective embeddings.
Indeed, the action on the normalization $\bP^1=\tilde{C}$ of $C$
is given by
$$
\varphi \, :\,    t\mapsto \frac{t}{1+at}.
$$
Note that 
the underlying topological space of $C$ is just $\bP^1$ and that 
all local rings coincide, {\em except} at the cusp $0$. 
The local ring 
${\cal O}_{C,0}$ is equal to the ideal generated by
$t^2,t^3$ in the ring $F[t]|_{t=0}$. 
This ideal is fixed under the action 
of $\varphi$. Therefore, the action descends to $C$.
By (\ref{prop:p1-all}), $C$ does not
admit a $\bG_a^1$-equivariant projective embedding.
\item 
Describe versal deformation spaces of 
non-normal proper $\bG_a^1$-varieties together with 
the $\bG_a^1$-action on these spaces.
\end{enumerate}

\section{Surfaces}

Throughout this section $X$ will be a 
smooth proper $\bG_a^2$-variety.

\begin{prop}\label{prop:contr}
Let $E\subset X$ be a $(-1)$-curve. 
Then there exists a morphism of
$\bG_a^2$-varieties $X\ra X'$ which blows down $E$.
\end{prop}

{\em Proof.} This follows from proposition 
\ref{prop:linearization} and corollary \ref{coro:linear}.

\begin{prop} \label{prop:hirz}
Every $\bG_a^2$-surface $X$ admits a $\bG_a^2$-equivariant
morphism onto $\bP^2$ or a Hirzebruch surface ${\bF}_n$. 
\end{prop}

{\em Proof.} This follows from the existence of minimal models  
for rational surfaces.

\subsection{Hirzebruch surfaces}
 
Let $X$ be a ${\bG}^2_a$-variety. Assume that $X$ is 
isomorphic
to ${\bF}_n$ as an algebraic variety with $n>0$.    
Its zero-section $e$ is stabilized under the group
action, as is the distinguished fiber $f$.
Let $\xi_n$ be the line bundle on $X$ corresponding to the 
section at infinity.  There is an induced equivariant morphism
$$\mu: X \rightarrow {\bP}(H^0(X,\xi_n)^*)={\bP}^{n+1}.$$
The image $\mu(X)$ is the cone over a smooth rational normal curve
of degree $n$;  $\mu$ contracts $e$ to the vertex of this cone.
 
We compute the representation of ${\bG}^2_a$ on $H^0(X,\xi_n)^*$.  
It has a distinguished one-dimensional fixed subspace
$W_1$ corresponding to the vertex.  The resulting
representation on $H^0(X,\xi_n)^*/W_1$  can be
easily understood geometrically.  
It has a one dimensional kernel
and the corresponding faithful representation
$${\bG}^1_a \rightarrow \GL (H^0(X,\xi_n)^*/W_1)$$
is the $n$-fold symmetric power of the standard 
two-dimensional representation.  
Here we are using the fact that $\mu(X)$
is the cone over a rational normal curve of degree $n$.
 
Choose a basis $S_1,S_2 \in L({\bG}^2_a)$ 
such that $S_1$ acts nontrivially 
and $S_2$ acts trivially as matrices on $H^0(X,\xi_n)^*/W_1$.
As a matrix on $H^0(X,\xi_n)^*$, 
$S_2$ has image  contained in $W_1$, 
and $S_1S_2=0$.  We have already seen that 
$S_1^n\ne 0$ as a matrix on $H^0(X,\xi_n)^*/W_1$.  

We consider two possible cases:  either $S_1^{n+1} \ne 0$
or $S_1^{n+1}=0$.  In the first case, we apply the following
fact about nilpotent matrices.
 
\begin{lem}
Let $S_1$ be an $(n+2)\times (n+2)$ nilpotent matrix such that
$S_1^{n+1} \ne 0$.  Then the centralizer of $S_1$ consists
of the algebra of matrices generated by $S_1$ and the identity.  
\end{lem}
 
{\em Proof.} This follows from a straightforward induction 
once we put $S_1$ in Jordan canonical form.
 
The lemma implies that $S_2$ may 
be written as some polynomial of $S_1$.  
The fact that the images of  $S_2$ and  $S^{n+1}_1$ 
both lie in  $W_1$ implies that $T=cS^{n+1}$ for 
some $c\ne 0$.  In this case
we have
$$H^0(X,\xi_n)^*=\rho_R \text{ where } 
R=F[S_1,S_2]/[S_1S_2,S_2-S_1^{n+1}].$$
Furthermore, $S_2$ acts nontrivially on the distinguished fiber $f$.

We now assume that $S_1^{n+1}=0$.  In this case we have
$$H^0(X,\xi_n)^*=\rho_R \text{ where } 
R=F[S_1,S_2]/[S_1S_2,S_1^{n+1}],$$
and the action is trivial along the distinguished fiber.  

In conclusion:
\begin{prop}
Let $X$ be a $\bG^2_a$-variety as above and let
$\xi_n$ denote the line bundle corresponding
to the section at infinity.  If the action on the distinguished
fiber is nontrivial then the representation
$${\bG}^2_a \rightarrow \GL (H^0(X,\xi_n)^*)=\GL_{n+2}$$
is equivalent to
$\exp (a_1 S_1+ a_{n+1}S_2)$ 
where $S_2=S_1^{n+1}\ne 0$ and $S_1S_2=0$.  
If the distinguished fiber is fixed under the action
then the representation is equivalent to
$\exp (a_1 S_1+ a_{n+1}S_2)$ where $S_1^{n+1}=0$ and $S_1S_2=0$.    
\end{prop}
 
A geometrical interpretation is obtained as follows.  
Let $W$ be the 
$(n+2)$-dimensional representation of ${\bG}^2_a$ described above.
Then each $X$ admits an equivariant birational
morphism into ${\bP}(W)$ and corresponds to the closure
of some nondegenerate orbit.  
To classify the surfaces $X$ it suffices to classify the 
nondegenerate ${\bG}_a^2$ orbits of 
${\bP}(W)$ modulo automorphisms, i.e.
the ${\bG}^2_a$-automorphisms of ${\bP}(W)$.

In the first case, these are exactly the automorphisms commuting
with the action of $S_1$,
i.e., the homotheties and the matrices
$$\exp (a_1 S_1 + a_2 S_1^2 + \ldots + a_{n+1}S_1^{n+1}).$$
Note that this gives ${\bP}(W)$ the structure of a 
${\bG}_a^{n+1}$-variety, which has a dense open orbit
equal to the complement of the distinguished hyperplane.  
In particular, any two nondegenerate orbit closures in ${\bP}(W)$
are related by an automorphism of ${\bP}(W)$.
It follows that $X$ is unique up to equivalence.  
 
In the second case, these automorphisms include the homotheties
and the matrices
$$\exp (a_1 S_1 + a_2 S_1^2 + \ldots + a_nS_1^n+a_{n+1}S_2).$$
Again, ${\bP}(W)$ has the structure of a ${\bG}_a^{n+1}$-variety,
and any two nondegenerate orbit closures are related by an automorphism.
It follows that $X$ is unique up to equivalence.  

These arguments yield the following: 
\begin{prop}
The Hirzebruch surfaces 
${\bF}_n$ with $n>0$ each have two distinct ${\bG}^2_a$-structures, 
including a unique structure with a 
nontrivial action on the distinguished
fiber.  The second structure is obtained by taking an elementary  
transformation of the structure on ${\bF}_{n-1}$.
The product 
${\bF}_0={\bP}^1 \times {\bP}^1$ has a unique ${\bG}^2_a$-structure,
induced from the ${\bG}^1_a$ actions on each factor.  
\end{prop}
The elementary transformation involves blowing-up the intersection
of the zero section and the distinguished fiber of ${\bF}_{n-1}$,
and then blowing-down the proper transform of this fiber. 
To prove the last statement, we project from
a fixed point of ${\bF}_0$.  The image is ${\bP}^2$ 
with ${\bG}_a^2$ acting by translation.  

\subsection{Examples and questions}

\begin{enumerate}
\item Interesting {\em singular} surfaces
admitting  $\bG_a^2$-structures: Del Pezzo 
surface of degree 5 with an isolated $A_4$-singularity. 
\item Can the $\bG_a^2$-structures on a given (smooth) 
surface have moduli? 
\item Classify $\bG_a^2$-structures on projective surfaces
with log-terminal singularities and Picard number 1. 
\end{enumerate}

\section{Threefolds}

\begin{theo}\label{theo:3-folds}
Let  $X$ be a smooth projective  $\bG_a^3$-variety
such that the boundary $D$ is irreducible.  Then $X$ 
is one of the following:
\begin{enumerate}
\item $X=\bP^3$, $D$ a hyperplane (the possible $\bG_a^3$-structures
were listed in \ref{prop:p3});
\item $Q_3\subset \bP^4$ is a smooth quadric, $D$ a tangent hyperplane section.
It has a unique $\bG_a^3$-structure (described in the proof). 
\end{enumerate}
\end{theo}

{\em Proof.} 
We know that $-K_X=r\cdot D$ where 
$r\ge 2$ (\ref{theo:anticanonical}).
Therefore, $X$ is a Fano variety of index $r\ge 2$ and it is rational.
Furthermore, if it has index $=2$ then the subgroup 
fixing $D$
has dimension one (by Corollary\ref{coro:Qncrit}).  We first consider
the case where the index $>2$.  

We show there is a unique ${\bG}^3_a$-structure
on $Q_3$, and that the boundary $D$ is necessarily
equal to a tangent hyperplane section.  
First let us convince ourselves that a quadric with a tangent
hyperplane is equivariant. 
Consider $\bP^3$ with the translation action $\tau_3$,
which fixes the hyperplane at infinity $P$.
Blow up a smooth conic curve  $C\subset P\subset \bP^3$ 
and blow down the proper transform of $P$. 
Now we prove that there are no other
$\bG_a^3$-structures on $Q_3$. Any $\bG_a^3$-action
on $Q_3$ has a fixed point $p$. Projecting 
from $p$ gives an equivariant 
birational map $f:Q_3 \dashrightarrow  \bP^3$.  
The induced map from the blow-up of $Q_3$ at $p$ to $\bP^3$ is
the blow-up of a conic contained in the proper transform of $D$
(where $D$ is the boundary of $\bG_a^3$ in $Q_3$).
The proper transform of $D$ is a plane $P$. 
The classification of $\bG_a^3$-structures on $\bP^3$ implies
that $P$ is {\em fixed} under 
the action of $\bG_a^3$ on $\bP^3$ (cf. \ref{prop:p3}).  

We return to the case where index equals two.  
By Furushima's classification of non-equivariant compactifications
of $\bG_a^3$ (cf. \cite{furushima}, \cite{furushima-2}), $X$ is 
a codimension 3 linear section of the Grassmannian ${\rm Gr}(2,5)$.
Consider the action induced on $F$, the variety of lines on $X$.
Any line on $X$ has normal bundle equal to ${\cO} \oplus {\cO}$
or ${\cO}(-1) \oplus {\cO}(+1)$ and there is always a line
of the second type (cf. \cite{fn}).  Choose such a line $L$ stable under
the ${\bG}_a^3$ action (we are choosing a fixed point on the locus of
lines of the second type in $F$).  
Projecting from $L$ gives an equivariant birational map
$$\pi_L: X \dashrightarrow Q_3$$
(cf. pp. 112 of \cite{fn}).

By \cite{furushima-2}, 
there are two cases to consider. 
In the first case the boundary $D\subset X$ is non-normal, with 
singular locus $L$. 
The total transform of $D$ consists of a
hyperplane section $H\subset Q_3$. 
The image of the boundary $D$ is a
smooth rational curve of degree 3, contained in $H$. 
We have already shown that there is a unique $\bG_a^3$-structure
on $Q_3$, which does not admit any smooth rational
curves of degree 3 contained in the boundary and stable
under the action (the only stable curve in $H$
is the distinguished ruling). 

In the second case the boundary 
$D\subset X$ is a normal singular Del Pezzo surface of degree 5
with an isolated $A_4$-singularity.  The curve $L\subset D$
is the unique $(-1)$-curve in the minimal resolution of $D$.  
Under $\pi_L$, $D$ is mapped birationally (and equivariantly!) 
to a tangent hyperplane
section $H \subset Q_3$.  The subgroup fixing $H$
has dimension 2, so the same
holds for $D$.  This contradicts Corollary \ref{coro:Qncrit}.

\subsection{Examples and questions}

\begin{enumerate}
\item
A singular projective $\bG_a^3$-variety with
one irreducible boundary divisor on which the $\bG_a^3$-action is
trivial. Construction: Blow up a pair 
of intersecting lines in $\bP^3$. Then 
blow down the proper transform of 
the plane containing them. The resulting
variety is a singular quadric hypersurface in $\bP^4$. 
\item  $\bG_a^3$-equivariant flop: Consider 
a quadric hypersurface 
$Q^*\subset \bP^4$ with an isolated singularity $p$.
Let $Y$ be the blow up of $Q^* $ at $p$. 
Let $E\equiv \bP^1\times \bP^1$ be the exceptional divisor.
Blowing down $E$ in different directions yields smooth threefolds
$X_1$ and $X_2$. It suffices to exhibit a $\bG_a^3$-structure
on $Q^*$ (all the constructions are natural and equivariant). 
This structure is obtained by using  
the fact that  $Q^*$ is a cone over $\bF_0$. 
\item $\bG_a^3$-equivariant flip:
Consider $\bF_1$ embedded in $\bP^4$ as a cubic scroll.
Let $X^*$ denote the cone over it. 
It is easy to see that $X^*$ is a $\bG_a^3$-variety.
Let $X$ be the blow-up
of $X^*$ at the singular point.  
The exceptional divisor $E$ is isomorphic to $\bF_1$. 
We can blow down $E$ along its ruling, 
which gives a non-singular threefold $Y$. 
The image of $E$ in $Y$ is a rational curve $C$ which
can be flipped $\bG_a^3$-equivariantly. 
\item 
If $X$ is a smooth projective $\bG_a^3$-variety and 
the action on the boundary is trivial is $X=\bP^n$? 
\end{enumerate}

\noindent Department of Mathematics \\
University of Chicago \\
5734 University Avenue \\
Chicago, IL 60637 \\
USA \\
{\em email address:} hassett@math.uchicago.edu

\

\

\noindent Department of Mathematics \\
U.I.C. \\
851 S. Morgan Str.\\
Chicago, IL 60607 \\
USA \\
{\em email address:} yuri@math.uic.edu

\end{document}